\documentclass[11pt, twoside]{amsart}

\usepackage{t1enc}
\usepackage[latin1]{inputenc}
\usepackage{latexsym}
\usepackage{amssymb}
\usepackage{graphicx}
\usepackage{amsmath}
\usepackage{amsthm}
\usepackage{amsfonts}
\usepackage{mathpazo}
\usepackage{mathrsfs}
\usepackage[paper = letterpaper, left = 3.6cm, right = 3.6cm, headsep = 6mm, footskip = 10mm,
top = 34mm, bottom = 36mm, footnotesep=5mm, headheight = 2cm]{geometry}
\usepackage{fancyhdr}
\usepackage[all]{xy}
\usepackage[british]{babel}
\usepackage{soul}
\usepackage{hyperref}
\usepackage[active]{srcltx}

\newcommand*{\widebar}{\overline}
\newcommand*{\definiere}{\mathrel{\mathop:}=}

\newcommand*{\tensor}{\otimes}

\newcommand{\cyclic}{\mathop{\kern0.9ex{{+}\kern-2.10ex\raise-0.20
      ex\hbox{\Large\hbox{$\circlearrowright$}}}}\limits}
\newcommand{\acts}{\mbox{ \raisebox{0.26ex}{\tiny{$\bullet$}} }}

\def\N{\ifmmode{\mathbb N}\else{$\mathbb N$}\fi}
\def\Z{\ifmmode{\mathbb Z}\else{$\mathbb Z$}\fi}
\def\Q{\ifmmode{\mathbb Q}\else{$\mathbb Q$}\fi}
\def\R{\ifmmode{\mathbb R}\else{$\mathbb R$}\fi}
\def\C{\ifmmode{\mathbb C}\else{$\mathbb C$}\fi}
\def\K{\ifmmode{\mathbb K}\else{$\mathbb K$}\fi}
\def\P{\ifmmode{\mathbb P}\else{$\mathbb P$}\fi}
\def\g{\ifmmode{\mathfrak g}\else {$\mathfrak g$}\fi}
\def\h{\ifmmode{\mathfrak h}\else {$\mathfrak h$}\fi}
\def\a{\ifmmode{\mathfrak a}\else {$\mathfrak a$}\fi}
\def\k{\ifmmode{\mathfrak k}\else {$\mathfrak k$}\fi}
\def\p{\ifmmode{\mathfrak p}\else {$\mathfrak p$}\fi}
\def\b{\ifmmode{\mathfrak b}\else {$\mathfrak b$}\fi}
\def\n{\ifmmode{\mathfrak n}\else {$\mathfrak n$}\fi}
\def\m{\ifmmode{\mathfrak m}\else {$\mathfrak m$}\fi}
\def\t{\ifmmode{\mathfrak t}\else {$\mathfrak t$}\fi}
\def\O{\ifmmode{\mathscr{O}}\else {$\mathscr{O}$}\fi}
\def\W{\ifmmode{\mathcal{V}}\else {$\mathscr{W}$}\fi}

\def\hq{/\hspace{-0.14cm}/}
\def\kleinematrix#1,#2,#3,#4,{\begin{pmatrix}#1 & #2 \\ #3 & #4
  \end{pmatrix}}

\DeclareMathOperator{\Lie}{Lie}

\DeclareMathOperator{\supp}{supp}

\DeclareMathOperator{\Ext}{Ext}
\DeclareMathOperator{\Hom}{Hom}
\DeclareMathOperator{\depth}{depth}

\newtheoremstyle{daniel}{3.0mm}{0mm}{\itshape}{}{\bfseries}{.}{1.5mm}{}
\theoremstyle{daniel}
\newtheorem{thm}{Theorem}[section]
\newtheorem{prop}[thm]{Proposition}
\newtheorem{Defi}[thm]{Definition}
\newtheorem{lemma}[thm]{Lemma}

\newtheorem{Exs}[thm]{Examples}
\newtheorem{Ex}[thm]{Example}
\newtheorem{Rems}[thm]{Remarks}
\newtheorem{Rem}[thm]{Remark}

\newtheorem*{thm*}{Theorem}
\newtheorem*{cor*}{Corollary}
\newtheorem*{thm3.6}{Theorem 3.6}
\newtheorem*{thm4.3}{Theorem 4.3}
\newtheorem*{prop*}{Proposition}

\newenvironment{rem}   {\begin{Rem}\em}{\end{Rem}}

\newenvironment{ex}  {\begin{Ex}\em}{\end{Ex}}
\newenvironment{exs}  {\begin{Exs}\em}{\end{Exs}}

\numberwithin{equation}{section}
 \def\polhk#1{\setbox0=\hbox{#1}{\ooalign{\hidewidth
  \lower1.5ex\hbox{`}\hidewidth\crcr\unhbox0}}}
  \def\polhk#1{\setbox0=\hbox{#1}{\ooalign{\hidewidth
  \lower1.5ex\hbox{`}\hidewidth\crcr\unhbox0}}}
\providecommand{\bysame}{\leavevmode\hbox to3em{\hrulefill}\thinspace}
\providecommand{\MR}{\relax\ifhmode\unskip\space\fi MR }

\providecommand{\href}[2]{#2}
\begin{document}
\title[Rational singularities and quotients]{Rational singularities and quotients by\\ holomorphic group actions}
\author{Daniel Greb}
\address{Albert-Ludwigs-Universit\"at\\
Mathematisches Institut\\
Abteilung f\"ur Reine Mathematik\\
Eckerstr.~1\\
79104 Freiburg im Breisgau\\
Germany}
\email{daniel.greb@math.uni-freiburg.de}{}
\urladdr{\href{http://home.mathematik.uni-freiburg.de/dgreb}{http://home.mathematik.uni-freiburg.de/dgreb}}

\thanks{\emph{Mathematical Subject Classification:} 32M05, 32S05, 32C36, 14L30}
\thanks{\emph{Keywords:} group actions on complex spaces, rational singularities, vanishing theorems}
\date{}

\begin{abstract}
We prove that rational and $1$-rational singularities of complex spaces are stable under taking quotients by holomorphic actions of reductive and compact Lie groups. This extends a result of Boutot to the analytic category and yields a refinement of his result in the algebraic category. As one of the main technical tools vanishing theorems for cohomology groups with support on fibres of resolutions are proven.
\end{abstract}
\maketitle
\section{Introduction and statement of main results}
Generalising the Hochster-Roberts Theorem \cite{HochsterRoberts}, Jean-Fran\c{c}ois Boutot proved that the class of varieties with rational singularities is stable under taking good quotients by reductive groups, see \cite{Boutot}. This result has found various applications in algebraic geometry, symplectic geometry, and representation theory.

In this note we study singularities of semistable quotients of complex spaces by complex reductive and compact Lie groups and prove a refined analytic version of Boutot's result.

The notion of \emph{semistable quotient} or \emph{analytic Hilbert quotient} is a natural generalisation of the classical concept of \emph{good quotient}, as used in Geometric Invariant Theory, to the holomorphic setup. Analytic Hilbert quotients play an important role in the study of group actions of complex reductive and compact groups on Stein and Kählerian complex spaces, cf.\ \cite{HeinznerGIT} and \cite{ReductionOfHamiltonianSpaces}. Here we show that the following refinement of Boutot's result holds for these quotients.
\begin{thm}\label{thm:maintheorem}
Let $G$ be a complex reductive Lie group acting holomorphically on a complex space $X$. Assume that the analytic Hilbert quotient $\pi: X \to X\hq G$ exists.
\begin{enumerate}
\item If $X$ has rational singularities, then $X\hq G$ has rational singularities.

\item If $X$ has 1-rational singularities, then $X\hq G$ has 1-rational singularities.
\end{enumerate}
\end{thm}
The notion of \emph{1-rational singularity} gains its importance from the fact that it provides the natural setup in which projectivity results for K\"ahler Moishezon manifolds extend to singular spaces, see \cite[Thm.\ 1.6]{ProjectivityofMoishezon}. The fundamental properties of $1$-rational singularities that allow to carry over the projectivity results from the smooth to the singular case were described in \cite[§12.1]{KollarMoriFlips}. Complex spaces with $1$-rational singularities also play a decisive role in the author's algebraicity results for momentum map quotients of algebraic varieties, see \cite{PaHq}.

As a corollary of Theorem~\ref{thm:maintheorem} we prove the following result for spaces on which there are a priori no actions of complex groups, e.g.\ bounded domains.
\begin{thm}\label{thm:mainthmK}
Let $K$ be a compact real Lie group acting holomorphically on a complex space $X$. Assume that the analytic Hilbert quotient $\pi: X \to X\hq K$ exists.
\begin{enumerate}
\item If $X$ has rational singularities, then $X\hq K$ has rational singularities.

\item If $X$ has 1-rational singularities, then $X\hq K$ has 1-rational singularities.
\end{enumerate}
\end{thm}
Furthermore,  we obtain the following refinement of Boutot's result in the algebraic category as a byproduct of our proof of Theorem~\ref{thm:maintheorem}.
\begin{thm4.3}
Let $X$ be a normal algebraic $G$-variety with good quotient $\pi: X \to X\hq G$. Let $f: \widetilde X \to X$ be a resolution of $X$, and let $g: Z \to X\hq G$ a resolution of $X\hq G$. Assume that $R^jf_*\mathscr{O}_{\widetilde X} = 0$ for $1 \leq j \leq q$. Then also $R^jg_*\mathscr{O}_Z = 0$ for $1 \leq j \leq q$.
\end{thm4.3}
This theorem generalises the main result of \cite{GrebSingularities} and can be used to study quotients of $G$-varieties with worse than rational singularities; see for example \cite[Lemma~3.3]{KovacsDuBoisI} for cohomological conditions imposed by Cohen-Macaulay singularities.

Although the main result refers to complex spaces with group action, one of the key technical results of this paper does not refer to the equivariant setup at hand and may also be of independent interest.
\begin{thm3.6}[Vanishing for cohomology with support I]
Let $f: \widetilde{X} \to X$ be a resolution of an irreducible normal complex space $X$ of dimension $n \geq 2$. Let $x \in X$, $F=f^{-1}(x)_{\mathrm{red}}$, and assume that $\bigcup_{k \geq 1} \supp R^kf_*\mathscr{O}_{\widetilde{X}} \subset \{x\}$. Then, we have
\[H^{q}_{F} \bigl(\widetilde{X}, \mathscr{O}_{\widetilde{X}} \bigr) = \{0\} \quad \text{ for all } q < n.\]
\end{thm3.6}
This vanishing theorem and the detailed analysis of algebraic and analytic cohomology groups with support generalises work of Karras~\cite{KarrasLocalCohomology}, cf. Section~\ref{sect:cohomologyvanishing}.

The author gratefully acknowledges the financial support of the Mathematical
Sciences Research Institute, Berkeley, during the 2009
"Algebraic Geometry" program. He wants to thank S\'andor Kov\'acs and Miles Reid for useful and stimulating discussions.
\section{Setup and notation}\label{sect:setup}
In the following a \emph{complex space} is a reduced complex space with countable topology. The structure sheaf of a given complex space $X$ will be denoted by $\mathscr{O}_X$. Furthermore, \emph{analytic subsets} of complex spaces are assumed to be closed. Throughout the paper $G=K^\C$ will always denote a complex reductive Lie group with maximal compact subgroup $K$.
\subsection{Singularities and their resolutions}\label{sect:sing}
Let $X$ be a complex space. Then, a \emph{resolution} of $X$ is a proper modification $f: \widetilde X \to X$ of $X$ such that $\widetilde X$ is a complex manifold.

If $f: \widetilde X \to X$ is a resolution of $X$, then for every $q$ the $q$-th push forward sheaf $R^qf_*\mathscr{O}_{\widetilde X}$ is a coherent sheaf of $\mathscr{O}_X$-modules. A complex space is said to have \emph{rational singularities} if it is normal and if for every resolution $f: \widetilde X \to X$ of $X$, the direct image sheaves $R^qf_*\mathscr{O}_{\widetilde X }$ vanish for $q > 0$. Analogously, we say that $X$ has \emph{1-rational singularities} if $X$ is normal, and if for every resolution $f: \widetilde X \to X$ the first direct image sheaf $R^1f_*\mathscr{O}_{\widetilde X}$ vanishes.

If $f_1: X_1 \to X$ and $f_2: X_2 \to X$ are two resolutions of $X$, then for every $q$ the sheaves $R^q(f_1)_*\mathscr{O}_{X_1}$ and $R^q(f_2)_*\mathscr{O}_{X_2}$ are isomorphic. In particular, having ($1$-)rational singularities is a local property that can be checked on a single resolution. Furthermore, if a Lie group $G$ acts on $X$ by holomorphic transformations, and if $f: \widetilde X \to X$ is a resolution, then the support of $R^qf_*\mathscr{O}_{\widetilde X } $ is a $G$-invariant analytic subset of $X$.

\begin{rem}\label{rem:rationalCM}
A complex space $X$ has rational singularities if and only if the following two conditions are satisfied:
\begin{enumerate}
\item $X$ is Cohen-Macaulay,
\item if $f: \widetilde X \to X$ is any resolution, then $f_*\omega_{\widetilde X} = \omega_X$.
\end{enumerate}
Here, $\omega_X, \omega_{\widetilde X}$ denote the \emph{canonical sheaves} of $X$, $\widetilde X$, respectively: if $Z$ is any normal complex space and if $\imath: Z_{\mathrm{smooth}} \to Z$ denotes the inclusion of the smooth part of $Z$ into $Z$, then $\omega_Z \definiere \imath_*(\Omega_{Z_{\mathrm{smooth}}}^{\dim Z})$.
\end{rem}

Since $R^{\dim X}f_*\mathscr{O}_{\widetilde X}$ vanishes for every resolution $f: \widetilde X \to X$ of a (normal) complex space $X$, the two notions rational and $1$-rational singularity coincide for two-di\-men\-sio\-nal complex spaces. In higher dimensions these notions differ as the following example shows.

\begin{ex}
Let $X$ be the affine cone over a smooth quartic hypersurface $S$ in $\P_3$. Then, the vertex of the cone is an isolated $1$-rational singularity that is not rational. In fact, if $f: \widetilde X \to X$ is the blowdown of the zero section in the total space $\widetilde X$ of the line bundle associated to $\mathscr{O}_S(-1)$, we have $(R^2f_*\mathscr{O}_{\widetilde X})_{\mathrm{vertex}} = \C$.
\end{ex}

\subsection{Analytic Hilbert quotients}
If $G$ is a real Lie group, then a \emph{complex $G$-space $Z$} is a
complex space with a real-analytic action $\alpha: G \times Z \to Z$ such that all the maps
$\alpha_g: Z \to Z$, $z \mapsto \alpha(g,z) = g \acts z$ are holomorphic. If $G$ is a complex Lie group, a
\emph{holomorphic $G$-space
$Z$} is a complex $G$-space such that the action map $\alpha: G \times Z \to Z$ is
holomorphic.

Let $G$ be a complex reductive Lie group and $X$ a holomorphic $G$-space. A
complex space $Y$ together with a $G$-invariant surjective holomorphic map $\pi: X \to Y$ is called
an \emph{analytic Hilbert quotient} (or \emph{semistable quotient} in the terminology of \cite{SemistableQuotients}) of $X$ by the action of $G$ if
\vspace{-2.5mm}
\begin{enumerate}
 \item $\pi$ is a locally Stein map, and
 \item $(\pi_*\mathscr{O}_X)^G = \mathscr{O}_Y$ holds.
\end{enumerate}
\vspace{-1mm}
Here, \emph{locally Stein} means that there exists an open covering of $Y$ by open Stein subspaces
$U_\alpha$ such that $\pi^{-1}(U_\alpha)$ is a Stein subspace of $X$ for all $\alpha$. By
$(\pi_*\mathscr{O}_X)^G$ we denote the sheaf $U \mapsto \mathscr{O}_X(\pi^{-1}(U))^G = \{f \in
\mathscr{O}_X(\pi^{-1}(U)) \mid f \;\; G\text{-invariant}\}$, $U$ open in
$Y$.

An analytic Hilbert quotient of a holomorphic $G$-space $X$ is unique up to biholomorphism once it exists and we will denote it by $X\hq G$. It has the following properties (see \cite{SemistableQuotients}):
\begin{enumerate}
\item Given a $G$-invariant holomorphic map $\phi: X \to Z$ into a
complex
space $Z$, there exists a unique holomorphic map $\widebar \phi: Y \to Z$ such that
$\phi = \widebar \phi \circ \pi$.

\item For every Stein subspace $A$ of $X\hq G$ the inverse image $\pi^{-1}(A)$ is a Stein subspace
of $X$.
\item If $A_1$ and $A_2$ are $G$-invariant analytic subsets of $X$, then
$\pi(A_1) \cap \pi(A_2) = \pi(A_1 \cap A_2)$.
\end{enumerate}

\begin{exs}
1.) If $X$ is a holomorphic Stein $G$-space, then the analytic Hilbert quotient exists and has the properties listed above (see \cite{HeinznerGIT} and \cite{Snow}).

2.) Let $X$ be a Kählerian complex space on which $G=K^\C$ acts holomorphically. Assume that the induced $K$-action is Hamiltonian with momentum map $\mu: X \to \Lie(K)^*$. If $X(\mu) = \{x \in X \mid \overline{G\acts x}\cap \mu^{-1}(0) \neq \emptyset \}$ is non-empty, it is $G$-invariant and open in $X$, and the analytic Hilbert quotient $\pi: X(\mu) \to X(\mu)\hq G$ exists, see \cite{ReductionOfHamiltonianSpaces}.
\end{exs}
\subsection{Quotients by compact Lie groups}
If $G=K^\C$ acts holomorphically on a Stein space $X$ with analytic Hilbert quotient $\pi: X\to X\hq G$, then the quotient $X\hq G$ can also be constructed without refering to the $G$-action: let $\mathscr{O}_X(X)^K = \mathscr{O}_X(X)^G$ denote the algebra of $K$-invariant holomorphic functions on $X$; associated with $\mathscr{O}_X(X)^K$ is the equivalence relation $\sim : \{(x_1, x_2) \in X \times X \mid f(x_1) = f(x_2)  \text{ for all } f \in \mathscr{O}_X(X)^K\}$ and the quotient $X/\negthickspace\sim$ with respect to this equivalence relation is  a complex space naturally isomorphic to $X\hq G$.

Therefore, one is lead to the following generalisation of the notion of analytic Hilbert quotient which allows to study holomorphic actions of compact groups $K$ that do not necessarily extend to actions of the complexification $K^\C$: let $K$ be a compact Lie group and let $X$ be a complex $K$-space; a complex space $Q$ together with a surjective holomorphic map $\pi: X \to Q$ is called an \emph{analytic Hilbert quotient} if every $q \in Q$ has an open Stein neighbourhood $U$ in $Q$ such that $\pi^{-1}(U)$ is Stein and such that the restriction $\pi|_{\pi^{-1}(U)}: \pi^{-1}(U) \to U$ induces an isomorphism of $\pi^{-1}(U)/\negthickspace\sim$ and $U$.

It is shown in \cite{HeinznerGIT} that for actions of compact groups on Stein spaces the analytic Hilbert quotient always exists and is a Stein space.
\section{Vanishing for cohomology with support}\label{sect:cohomologyvanishing}
The proof of the main result uses the vanishing of cohomology groups with support on a fibre of a resolution $\pi: \widetilde X \to X$, see Theorems~\ref{thm:GRvanishingI} and \ref{thm:GRvanishingII}. This vanishing will be derived from analytic Grauert-Riemenschneider vanishing using duality. It is the analytic analogue of a vanishing result already used by Boutot in his proof of the rationality of quotient singularities. However, in contrast to the algebraic case, the analytic setup requires a careful study of the a priori different algebraic and analytic cohomology groups with support. This section is devoted to the proof of Theorems \ref{thm:GRvanishingI} and \ref{thm:GRvanishingII}
\subsection{Algebraic and analytic cohomology with support}
If $Y$ is a closed subspace of a topological space $X$ and if $\mathscr{F}$ is a sheaf of abelian groups on $X$ we denote the \emph{$q$-th (analytic) cohomology group with support on $Y$} by $H^q_Y\bigl( X, \mathscr{F}\bigr)$. Furthermore, we denote the \emph{$q$-th cohomology sheaf with support on $Y$} by $\mathscr{H}_Y^q(\mathscr{F})$. Suppose now that $Y$ is an analytic subset of a complex space $X$ and that $\mathscr{F}$ is a coherent analytic sheaf on $X$. Then, $\mathscr{H}_Y^q(\mathscr{F})$ is an analytic $\mathscr{O}_X$-sheaf with support contained in $Y$. It is however not coherent in general. Let $\mathscr{I}$ be an ideal sheaf defining $Y$. Then, we define
\[\Gamma_{[Y]}\big(X, \mathscr{F} \bigr) \definiere \lim_{\overset{\longrightarrow}{m}}\Hom_{\mathscr{O}_X} \bigl(\mathscr{O}_X/\mathscr{I}^m, \mathscr{F} \bigr), \quad  H^q_{[Y]}\bigl(X, \mathscr{F} \bigr) \definiere \lim_{\overset{\longrightarrow}{m}} \Ext^q\bigl(X; \mathscr{O}_X/\mathscr{I}^m, \mathscr{F} \bigr). \]We call $ H^q_{[Y]}\bigl(X, \mathscr{F} \bigr) $ the $q$-th \emph{algebraic cohomology group with support on $Y$}.
For each $m$ there exists a natural map $\Hom_{\mathscr{O}_X}\bigl(\mathscr{O}_X/\mathscr{I}^m, \mathscr{F} \bigr) \to \Gamma_Y\bigl(X, \mathscr{F} \bigr)$. For every $q \geq 0$ this induces a natural homomorphism
\[\rho_q:  H^q_{[Y]}\bigl(X, \mathscr{F} \bigr)  \to  H^q_{Y}\bigl(X, \mathscr{F} \bigr).\]
The homomorphism $\rho_0$ is an isomorphism. In the algebraic category this fact is used to prove that $\rho_q$ is bijective for all $q\geq 0$, see \cite[Thm.\ 2.8]{HartshorneLocalCohomology}. In the analytic category this is no longer true in general. However, we prove the following comparison result:
\begin{prop}[Comparison between algebraic and analytic local cohomology]\label{prop:comparisonI}
Let $f: \widetilde X \to X$ be a resolution of an irreducible normal complex space $X$ of dimension $n \geq 2$. Let  $x \in X$, $F=f^{-1}(x)_{\mathrm{red}}$, and assume that $\bigcup_{k \geq 1} \supp R^kf_*\mathscr{O}_{\widetilde X} \subset \{x\}$. Then, the natural homomorphisms
\[\rho_q: H^q_{[F]}\bigl(\widetilde X, \mathscr{O}_{\widetilde X} \bigr) \to H^q_{F}\bigl(\widetilde X, \mathscr{O}_{\widetilde X} \bigr) \]
are isomorphisms for $q< n$.
\end{prop}
\subsubsection{Coherence of cohomology sheaves and comparison}
In the proof of Proposition~\ref{prop:comparisonI} we are going to use the following comparison result.
\begin{prop}[see\ Sect.\ 1.3 of \cite{KarrasLocalCohomology} and Sect.\ 5.5 of \cite{TrautmannEndlichkeitssatz}]\label{prop:isomorphismcondition}
Let $Y$ be an analytic subset of a complex space $X$ and let $\mathscr{F}$ be a coherent analytic sheaf on $X$. Assume that the analytic sheaves $\mathscr{H}_Y^j(\mathscr{F})$ are coherent for $0\leq j \leq q$. Then, the map $\rho_j$ is an isomorphism for $0\leq j\leq q$. If $Y = \{x\}$ is a point, then additionally $\rho_{q+1}$ is injective.
\end{prop}
Therefore, one is led to study the coherence of the local cohomology sheaves $\mathscr{H}_Y^j(\mathscr{F})$. We will do this using the following coherence criterion.
\begin{prop}[Thm.\ 3.5 in \cite{SiuTrautmannGapSheaves}]\label{prop:coherencecondition}
Let $Y$ be an analytic subset of a complex space $X$ and let $\mathscr{F}$ be a coherent analytic sheaf on $X$. For every $m \in \N$ let $S_m(\mathscr{F}) \definiere \{x \in X \mid \depth_x \mathscr{F} \leq m\}$. Then, for every $q\geq 0$ the following conditions are equivalent:
\begin{enumerate}
\item The sheaves $\mathscr{H}_Y^j(\mathscr{F})$ are coherent for $0\leq j \leq q$.

\item $\dim_\C \bigl(Y \cap \overline{S_{k+q+1}(\mathscr{F}|_{X \setminus Y}) }^X\bigr)\leq k$ for all $k \in \Z$.
\end{enumerate}
\end{prop}
\subsubsection{Spectral sequences associated to a resolution}Let $f:\widetilde X \to X$ be a resolution of an irreducible normal complex space, $x\in X$, $F \definiere f^{-1}(x)_{\mathrm{red}}$, and let $\mathscr{F}$ be a coherent analytic sheaf on $\widetilde X$.  Analogous to \cite[Sect.\ 1.5]{KarrasLocalCohomology} we have
\[\Gamma_x\big(X, f_*\mathscr{F} \bigr) \cong \Gamma _F\bigl(\widetilde X, \mathscr{F}  \bigr)  \quad \text{and} \quad  \Gamma_{[x]}\big(X, f_*\mathscr{F} \bigr) \cong \Gamma _{[F]}\bigl(\widetilde X, \mathscr{F}  \bigr).\]
These isomorphisms induce the following commutative diagram of spectral sequences
\begin{equation*}
\begin{xymatrix}{
H^j_{[x]}\bigl(X, R^kf_*\mathscr{F} \bigr)\ar^{\rho_{jk}}[r]\ar@{=>}[d] & H^j_{x}\bigl(X, R^kf_*\mathscr{F} \bigr) \ar@{=>}[d]\\
H^m_{[F]}\bigl( \widetilde X, \mathscr{F}\bigr) \ar^{\rho_m}[r]& H^m_{F}\bigl( \widetilde X, \mathscr{F} \bigr),
}\end{xymatrix}
\end{equation*}
where $m = j+k$.
\begin{proof}[Proof of Proposition~\ref{prop:comparisonI}]
Using the commutative diagram of spectral sequences above and analysing the proof of \cite[Thm.\ 5.2.12]{Weibel}  we see that it suffices to show that
\begin{enumerate}
\item $\rho_{j0}$ is bijective for all $j < n$, injective for $j=n$, and
\item $\rho_{jk}$ is bijective for all $j > 0$, $k > 0$.
\end{enumerate}
Since $X$ is assumed to be normal, we have $f_*\mathscr{O}_{\widetilde X} = \mathscr{O}_X$ . Furthermore, again by assumption, the support of $R^kf_*\mathscr{O}_{\widetilde X}$ is contained in $\{x\}$ if $k > 0$. It follows that the complex space $X \setminus \{x\}$ has rational singularities. In particular, $X\setminus \{x\}$ is Cohen-Macaulay, i.e.\ $\depth_y \mathscr{O}_X = n$ for all $y \in X \setminus \{x\}$, cf.\ Remark~\ref{rem:rationalCM}. Combining these observations with Proposition~\ref{prop:coherencecondition} we see that  $\mathscr{H}_x^j(f_* \mathscr{O}_{\widetilde X}) = \mathscr{H}_x^j(\mathscr{O}_X)$ is coherent for all $j<n$.

If $k>0$, the fact that $\supp R^kf_*\mathscr{O}_{\widetilde X} \subset \{x\}$ implies $\mathscr{H}_x^j(R^kf_*\mathscr{O}_{\widetilde X}) = 0$ for all $j > 0$, cf.\ \cite[Sect.\ 0.4]{SiuTrautmannGapSheaves}.
Now we use Proposition~\ref{prop:isomorphismcondition} to conclude that both (1) and (2) hold.
\end{proof}
\subsection{Duality for algebraic cohomology with support}\label{sect:duality}
In this section we prove a duality theorem for algebraic cohomology with support. If $\mathscr{G}$ is a coherent analytic sheaf on a complex space $X$ we denote by $\mathscr{G}^*$ its dual sheaf.
\begin{prop}[Duality for algebraic cohomology with support]\label{prop:duality}
Let $f: \widetilde X \to X$ be a resolution of an irreducible normal complex space $X$ of dimension $n$. Let  $x \in X$, $F=f^{-1}(x)_{\mathrm{red}}$, and let $\mathscr{G}$ a locally free coherent analytic sheaf on $\widetilde X$. If $\supp R^qf_*\mathscr{G} \subset \{x\}$, then
\[R^qf_*\mathscr{G}  \overset{\mathsf{dual}}{\sim} H^{n-q}_{[F]}\bigl(\widetilde{X}, \omega_{\widetilde X}\tensor \mathscr{G}^{*}\bigr) \quad \text{ for all } 0\leq q \leq n. \]
\end{prop}
\begin{proof}(cf.\ \cite[Prop.\ 11.6.1]{KollarSingularitiesOfPairs} and \cite[Prop. 2.2]{HartshorneOgus})
Let $F_m$ denote the $m$-th infinitesimal neighbourhood of the fibre $F$, i.e., if $\mathscr{I}_F$ denotes the ideal sheaf of $F$ in $\mathscr{O}_X$, we have $\mathscr{O}_{F_m} = \mathscr{O}_{\widetilde{X}}/\mathscr{I}_{F}^m$. Let $\imath_m : F_m \to \widetilde{X}$ be the inclusion. Since the dimension of the cohomology group $H^{q+1}\bigl(F_m, \imath_m^*\mathscr{G}\bigr) = H^{q+1}(\widetilde{X}, \mathscr{O}_{F_m}\otimes \mathscr{G}\bigr)$ is finite by the Cartan-Serre Theorem, we may apply Theorem I and Theorem 6 of \cite{AndreottiDuality} to conclude that for every $m \in \N$
\begin{equation*}
H^q\bigl(\widetilde{X}, \mathscr{O}_{F_m}\otimes \mathscr{G} \bigr) \overset{\mathsf{dual}}{\sim} \Ext_c^{n-q}\bigl(\widetilde{X}; \mathscr{O}_{F_m}, \omega_{\widetilde{X}}\otimes \mathscr{G}^* \bigr).
\end{equation*}
Here, $\Ext_c^q\bigl(\widetilde{X}; \mathscr{O_{F_m}}, \cdot \bigr)$ is the $q$-th Ext-group with compact support, cf.\ \cite[§6]{AndreottiDuality}.
Since $\supp \mathscr{O}_{F_m} = F$ is compact, it follows that
\begin{equation}\label{duality}
H^q\bigl(\widetilde{X}, \mathscr{O}_{F_m}\otimes \mathscr{G} \bigr) \overset{\mathsf{dual}}{\sim} \Ext^{n-q}\bigl(\widetilde{X}; \mathscr{O}_{F_m}, \omega_{\widetilde{X}}\otimes \mathscr{G}^* \bigr).
\end{equation}
Since the completion of $R^qf_*\mathscr{G}$ at $x$ is finite-dimensional by assumption, Grauert's comparison theorem (see \cite[Hauptsatz IIa]{GrauertComparison}) implies that for large $m$
\[H^q\bigl(\widetilde{X}, \mathscr{O}_{F_m} \otimes \mathscr{G} \bigr) \cong (R^qf_*\mathscr{G})_x = R^qf_*\mathscr{G}.\]
On the other hand, the duality~\eqref{duality} implies that $\Ext^{n-q}\bigl(\widetilde{X}; \mathscr{O}_{F_m}, \omega_{\widetilde{X}} \otimes \mathscr{G}^*\bigr)$ is finite-dimensional. Consequently, for large $m$, we obtain
\[H^{n-q}_{[F]}\bigl(\widetilde{X}, \omega_{\widetilde{X}} \otimes \mathscr{G}^* \bigr) =  \Ext^{n-q}\bigl(\widetilde{X}; \mathscr{O}_{F_m}, \omega_{\widetilde{X}} \otimes \mathscr{G}^*\bigr) \overset{\mathsf{dual}}{\sim} H^q\bigl(\widetilde{X}, \mathscr{O}_{F_m}\otimes \mathscr{G} \bigr) = R^qf_*\mathscr{G}. \qedhere \]
\end{proof}
\subsection{Vanishing for analytic cohomology with support}
Using the results obtained in the previous two sections, we are now able to deduce the following vanishing theorems, the main technical results of this note.
\begin{thm}[Vanishing for cohomology with support I]\label{thm:GRvanishingI}
Let $f: \widetilde{X} \to X$ be a resolution of an irreducible normal complex space $X$ of dimension $n \geq 2$. Let $x \in X$, $F=f^{-1}(x)_{\mathrm{red}}$, and assume that $\bigcup_{k \geq 1} \supp R^kf_*\mathscr{O}_{\widetilde{X}} \subset \{x\}$. Then,  we have
\[H^{q}_{F} \bigl(\widetilde{X}, \mathscr{O}_{\widetilde{X}} \bigr) = \{0\} \quad \text{ for all } q < n.\]
\end{thm}
\begin{proof}
Since $\bigcup_{k \geq 1} \supp R^kf_*\mathscr{O}_{\widetilde X} \subset \{x\} $ by assumption, we use Proposition~\ref{prop:comparisonI}, Proposition~\ref{prop:duality}, and the Grauert-Rie\-men\-schnei\-der Theorem in its analytic version proven by Silva \cite[A.2.Lemma]{SilvaRelativeVanishing} to obtain
\[H^{q}_{F} \bigl(\widetilde{X}, \mathscr{O}_{\widetilde{X}} \bigr) \cong H^{q}_{[F]} \bigl(\widetilde{X}, \mathscr{O}_{\widetilde{X}} \bigr) \overset{\mathsf{dual}}{\sim} R^{n-q}f_*\omega_{\widetilde{X}} = 0 \quad \text{   for all } q < n. \qedhere\]
\end{proof}
To handle complex spaces with $1$-rational singularities, we prove the following result.
\begin{thm}[Vanishing for cohomology with support II]\label{thm:GRvanishingII}
Let $f:\widetilde{X} \to X$ be a resolution of an irreducible normal Stein space $X$ of dimension $n \geq 2$. Let $x \in X$, $F=f^{-1}(x)_{\mathrm{red}}$, and assume that $ \supp R^1 f_*\mathscr{O}_{\widetilde{X}} \subset \{x\}$. Then, we have
\[H^{1}_{F} \bigl(\widetilde{X}, \mathscr{O}_{\widetilde{X}} \bigr) = \{0\}\]
\end{thm}
\begin{proof}Consider the exact sequence for local cohomology (\cite[§ 0]{SiuTrautmannGapSheaves})
\[0 \to  H^0\bigl(\widetilde X, \mathscr{O}_{\widetilde X} \bigr)  \overset{\alpha}{\to}H^0\bigl(\widetilde X \setminus F, \mathscr{O}_{\widetilde X} \bigr) \to H^1_F\bigl(\widetilde X, \mathscr{O}_{\widetilde X} \bigr) \to  H^1\bigl(\widetilde X, \mathscr{O}_{\widetilde X} \bigr) \to \cdots.\]
Since $X$ is assumed to be normal, the map $\alpha$ is bijective. Consequently, $H^1_F\bigl(\widetilde X, \mathscr{O}_{\widetilde X} \bigr) $ injects into $ H^1\bigl(\widetilde X, \mathscr{O}_{\widetilde X} \bigr) = H^0 \bigl(X, R^1f_*\mathscr{O}_{\widetilde X} \bigr)$, and is therefore finite-dimensional by the assumption on the support of $R^1f_*\mathscr{O}_{\widetilde X}$. Applying \cite[Thm.\ 2.3]{HartshorneAmple} we see that the duality homomorphism
\[ H^{n-1}\bigl(F, \omega_{\widetilde X}|_{F}\bigr) \to H^1_F\bigl(\widetilde X, \mathscr{O}_{\widetilde X} \bigr)^* \] is surjective, where $\cdot |_{F}$ denotes set-theoretic restriction. Since $H^{n-1}\bigl(F, \omega_{\widetilde X}|_{F}\bigr)$ is isomorphic to the stalk of $\bigl(R^{n-1}f_* \omega_{\widetilde X}\bigr)$ at $x$, the claimed vanishing again follows from analytic Grauert-Riemenschneider vanishing \cite{SilvaRelativeVanishing}.
\end{proof}
\section{Proof of the main results}\label{sect:proof}
\subsection{Cutting by hyperplane sections}
The following lemma will be used to obtain information about the singularities of a complex space $X$ from information about the singularities of a general hyperplane section $H$ of $X$ and vice versa. Recall that a subset of a topological space $X$ is called \emph{fat} if it contains a countable intersecton of dense subsets of $X$. Every fat subset of a complete metric space $X$ is dense in $X$.
\begin{lemma}\label{lem:hyperplanecut}
Let $X$ be a normal Stein space, $\dim X \geq 2$, and let $f: \widetilde X \to X$ be a resolution. Let $\mathscr H \subset \mathscr{O}_X(X)$ be a finite-dimensional linear subspace containing the constants, $\dim_\C\mathscr{H}=N+1$, such that the image of $X$ under the associated map $\varphi_{\mathscr H}: X \to \P_N = \P(\C \oplus \C^N)_{[Z_0: Z_1: \dots:Z_N]}$ is closed in $\{Z_0 \neq 0\}$ and fulfills $\dim \varphi_{\mathscr H}(X) \geq 2$. Then, there exists a fat subset $U \subset \P(\mathscr H)$ such that for all $h \in U$ the following holds for the corresponding hyperplane section $\{h = 0\} \subset X$:
\begin{enumerate}
\item The preimage $\widetilde {H} \definiere f^{-1}(H)$ is smooth and $f|_{\widetilde H}: \widetilde{H} \to H$ is a resolution of $H$.
\item We have $R^qf_*\mathscr{O}_{\widetilde{H}} \cong R^qf_*\mathscr{O}_{\widetilde X} \tensor \mathscr{O}_{\overset{}{H}}\,$ for all $q \geq 0$.
\end{enumerate}
\end{lemma}
\begin{proof}
(1) This follows from Bertini and Sard theorems, see e.g. \cite{ManaresiPermanence}.

(2) In the exact sequence
\begin{equation}\label{idealsequence}
0 \to \mathscr{O}_{\widetilde X}(-\widetilde H) \overset{m}{\to} \mathscr{O}_{\widetilde X} \to \mathscr{O}_{\widetilde H} \to 0,
\end{equation}
the map $m$ is given by multiplication with the equation $h\in \mathscr{O}_X(X)$ defining $H$ and $\widetilde{H}$. Pushing forward the short exact sequence \eqref{idealsequence} by $f_*$ yields the long exact seqence
\begin{align}\label{longsheafsequence}
\begin{split}
0 &\to f_* \mathscr{O}_{\widetilde X}(-\widetilde H) \overset{m_0}{\longrightarrow} f_*\mathscr{O}_{\widetilde X} \to f_*\mathscr{O}_{\widetilde H} \to \\
  &\to R^1f_* \mathscr{O}_{\widetilde X}(-\widetilde H) \overset{m_1}{\longrightarrow} R^1 f_*\mathscr{O}_{\widetilde X} \to R^1 f_*\mathscr{O}_{\widetilde H} \to \cdots \\
  &\to  R^qf_* \mathscr{O}_{\widetilde X}(-\widetilde H) \overset{m_q}{\longrightarrow} R^q f_*\mathscr{O}_{\widetilde X} \to R^q f_*\mathscr{O}_{\widetilde H} \to \cdots \\
  &\to R^{n-1}f_* \mathscr{O}_{\widetilde X}(-\widetilde H) \overset{m_{n-1}}{\longrightarrow} R^{n-1} f_*\mathscr{O}_{\widetilde X} \to 0.
\end{split}
\end{align}
The maps $m_q$, $q=0,1,2, \dots, n-1$, are given by multiplication with the element $h \in \mathscr{O}_X(X)$. Since $\mathscr{H}$ is base-point free, there exists a fat subset $U$ in $\P(\mathscr{H})$ such that all the maps $m_q$ in the sequence~\eqref{longsheafsequence} are injective, cf.~\cite{ManaresiPermanence}.

Since $\mathscr{O}_{\widetilde X}(- \widetilde{H}) \cong f^* (\mathscr{O}^{}_X(-H))$, the projection formula for locally free sheaves yields
\begin{align*}
R^qf_*\mathscr{O}_{\widetilde X}(- \widetilde H) \cong R^qf_*\mathscr{O}_{\widetilde X} \tensor \mathscr{O}_{X}^{}(-H).
\end{align*}
Furthermore, for every $q$, the image of $m_q$ coincides with the image $\mathscr{B}_q$ of the natural map $R^qf_*\mathscr{O}_{\widetilde X} \otimes \mathscr{O}_{X}^{}(-H) \to R^qf_*\mathscr{O}_{\widetilde X}$.  Since $m_q$ is injective, it follows that
\[R^qf_*\mathscr{O}_{\widetilde H} \cong R^qf_*\mathscr{O}_{\widetilde X} / \mathscr{B}_q \quad \text{for all } q=0,1, \dots, n-1.\]
Tensoring with $R^qf_*\mathscr{O}_{\widetilde X}$ is right-exact, and hence the exact sequence
$0 \to \mathscr{O}_X(-H) \to \mathscr{O}_X \to \mathscr{O}_H \to 0$
yields $R^q f_* \mathscr{O}_{\widetilde H} \cong R^qf_*\mathscr{O}_{\widetilde X} \tensor \mathscr{O}^{}_H$, as claimed.
\end{proof}
\subsection{Proof of Theorem~\ref{thm:maintheorem}}
Let $\pi: X \to X\hq G$ be an analytic Hilbert quotient. The following lemma relates cohomology groups of $X\hq G$ to cohomology groups on $X$. This is essentially the only point in the proof of Theorem~\ref{thm:maintheorem} at which we use the fact that $\pi$ is the analytic Hilbert quotient for the $G$-action on $X$.
\begin{lemma}\label{lem:pullbackinjective}
Let $X$ be a holomorphic $G$-space with analytic Hilbert quotient $\pi: X \to X\hq G$. Then, for all $q \geq 0$, the natural map
\[\pi^*: H^q\bigl(X\hq G, \mathscr{O}_{X/\negthickspace/ G} \bigr) \to H^q\bigl(X, \mathscr{O}_X \bigr)\]
is injective.
\end{lemma}
\begin{proof}

Let $ \{U_\alpha\}_{\alpha \in \mathbb A}$ be a \v{C}ech covering of $X\hq G$. Since $\pi$ is Stein, $\{\pi^{-1}(U_\alpha)\}_{\alpha \in \mathbb A}$ is a \v{C}ech covering of $X$. If $[(v_{i_0 \dots i_q})_{i_k \in \mathbb{A}}] \in \ker (\pi^*)$, then there exists a cocycle $(w_{j_o \dots j_{q-1}})_{j_k \in \mathbb{A}}$ with $(\pi^*(v_{i_0 \dots i_q} ))_{i_k \in \mathbb{A}} = \delta ((w_{j_o \dots j_{q-1}})_{j_k \in \mathbb{A}})=\delta\bigl((\int_K w_{j_o \dots j_{q-1}} \mathrm{dk} )_{j_k \in \mathbb{A}}\bigr)$. Here, $\int_K \mathrm{dk}$ denotes the integral with respect to the Haar measure on a maximal compact subgroup $K$ of $G$. Since $(\pi_* \mathscr{O}_X)^K= (\pi_* \mathscr{O}_X)^G = \mathscr{O}_{X/\negthickspace/G}$, this implies $[(v_{i_0 \dots i_q})_{i_k \in \mathbb{A}}]= 0 \in H^q\bigl(X\hq G, \mathscr{O}_{X/\negthickspace/G}\bigr)$.
\end{proof}
With the results of the previous sections at hand we are now in the position to carefully adapt Boutot's proof to our analytic situation. Furthermore, by thoroughly carrying out the reduction steps omited in \cite{Boutot} we will even find a refinement of the classical result in the algebraic category, cf.\ Theorem~\ref{cor:Boutotrefined}.
\begin{proof}[Proof of Theorem \ref{thm:maintheorem}]
Since the claim is local and $\pi$ is Stein, we may assume that $X\hq G$ and $X$ are
Stein, and that $X\hq G$ admits a closed embedding $\varphi: X\hq G \hookrightarrow \C^N$ for some $N \in \N$. Since $\mathscr{O}_{X/\negthickspace/G} = (\pi_*\mathscr{O}_X)^G$, normality of $X$ implies that $X\hq G$ is also normal. As a consequence, we can assume in the following that $X$ is $G$-irreducible.

We prove the claim by induction on $\dim X\hq G$. For $\dim X\hq G = 0$ there is nothing to show.
For $\dim X\hq G = 1$ we notice that $X\hq G$ is smooth. So, let $\dim X\hq G \geq 2$. Let
$\pi: X \to X\hq G$ denote the quotient map and let $p_X: \widetilde X \to X$ be a resolution of $X$.
First, we claim that there exists a fat subset of the set of all affine hyperplane sections $H \subset X\hq G$ such that $H$ has rational singularities. Indeed, by Lemma \ref{lem:hyperplanecut} there exists a fat subset $U$ of the set of all affine hyperplane sections $H \subset X\hq G$ such that $p_X|_{\widetilde H}: \widetilde H \to \pi^{-1}(H)$ is a resolution, where $\widetilde H = p_X^{-1}(\pi^{-1}(H))$. By the same lemma we may furthermore assume that for all $H\in U$ we have
\[R^q(p_X)_*\mathscr{O}_{\widetilde{X}} \otimes_{\mathscr{O}_X}\mathscr{O}_{\pi^{-1}(H)} =
R^q (p_X)_*\mathscr{O}_{\widetilde{H}} \quad \text{ for all } q\geq 0. \] Since $(p_X)_* \mathscr{O}_{\widetilde{X}}=\mathscr{O}_X$ by
the analytic version of Zariski's main theorem, it follows from the case $q=0$ that $\pi^{-1}(H)$ is normal. Alternatively, one could invoke Seidenberg's Theorem (see \cite{ManaresiPermanence}). Together with the cases $q=1, \dots, \dim X -1$ this implies that $\pi^{-1}(H)$ has rational singularities. By
induction, it follows that $H = \pi^{-1}(H)\hq G$ has rational singularities.

Let $p: Z \to X\hq G$ be a resolution of $X\hq G$. The previous considerations together with Lemma~\ref{lem:hyperplanecut} imply that there exists a fat subset $U$ of the set of all affine hyperplane sections $H$ of $X\hq G$ such that $H$ has rational singularities, such that the restriction of $p$ to $\widehat{H} \definiere p^{-1 }(H)$, $p|_{p^{-1}(H)}: \widehat H \to H$, is a resolution of $H$, and such that $\mathscr{O}_H \otimes R^qp_*\mathscr{O}_Z = R^qp_*\mathscr{O}_{\widehat{H}}=0$. Consequently, the support of
$R^qp_*\mathscr{O}_Z$ does not intersect any $H\in U$ and hence $\supp(R^qp_*\mathscr{O}_Z)$ consists of isolated points. Since the claim is local, we can therefore assume in the following that $\bigcup_{q\geq 1}R^qp_*\mathscr{O}_Z \subset \{x_0\}$ for some $x_0 \in X\hq G$.

The group $G$ acts on the fibre product $Z\times_{X /\negthickspace / G} X$ such that the canonical projection $p_X: Z \times_{X /\negthickspace/ G}
X \to X$ is equivariant. One of the $G$-irreducible components $\widetilde{X}$ of $Z\times_{X /\negthickspace/ G} X$
is bimeromorphic to $X$, and, by passing to a resolution of $\widetilde{X}$ if necessary, we can assume that
$p_X: \widetilde{X} \to X$ is a resolution of $X$. We obtain the following commutative diagram
\[\begin{xymatrix}{
X \ar[d]_\pi & \ar[l]_{p_X} \ar[d]^{p_Z}\widetilde{X} \\
X\hq G & \ar[l]_<<<<<{p} Z .
}
  \end{xymatrix}
\]
Since $\bigcup_{q\geq 1}\supp (R^qp_*\mathscr{O}_Z) \subset \{x_0\}$, we have $(R^qp_*\mathscr{O}_Z)_{x_0} = H^0(X\hq G, R^qp_*\mathscr{O}_Z)$ for every $q \geq 1$. Recall that $X\hq G$ is Stein, hence the Leray spectral sequence for $p$ implies that it suffices to show that $H^q\bigl(Z, \mathscr{O}_Z\bigr) = 0$ for all $q\geq 1$. Since $X$ is Stein and has rational singularities, it follows that $H^q\bigl(\widetilde{X}, \mathscr{O}_{\widetilde{X}}\bigr) \cong H^q\bigl(X, \mathscr{O}_X\bigr) = \{0\}$ for all $q\geq 1$. Consequently, it
suffices to show that there exists an injective map $H^q\bigl(Z, \mathscr{O}_Z\bigr) \hookrightarrow H^q\bigl(\widetilde{X}, \mathscr{O}_{\widetilde{X}}\bigr)$.

We introduce the following notation: $ U= (X\hq G) \setminus \{x_0\}$, $U' =\pi^{-1}(U) \subset X$, $
\widetilde{U}= p_X^{-1}(U') \subset \widetilde{X}$, $V = p^{-1}(U) \subset Z$.
For every $q \geq 1$ we obtain the following commutative diagram of canonical maps
\begin{equation}\label{diagramrational}
\begin{xymatrix}{
 H^q\bigl(\widetilde{X}, \mathscr{O}_{\widetilde{X}}\bigr) \ar[r]^{h_{\widetilde X , \widetilde U}^q}& H^q \bigl(\widetilde{U},
\mathscr{O}_{\widetilde{U}}\bigr) & \ar[l]_{h_{U',\widetilde U}^q} H^q\bigl(U', \mathscr{O}_{U'}\bigr) \\
H^q\bigl(Z, \mathscr{O}_Z\bigr) \ar[u]_{h_{Z,\widetilde{X}}^q}\ar[r]^{h_{Z,V}^q}& H^q\bigl(V, \mathscr{O}_V\bigr) \ar[u]&
H^q\bigl(U, \mathscr{O}_U\bigr)\ar[l]_{h_{U,V}^q}\ar[u]_{h_{U,U'}^q}.
}
  \end{xymatrix}
\end{equation}
Let $Y \definiere p^{-1}(x_0) \subset Z$. We consider the exact sequence for cohomology with support
\begin{equation}\label{localcohomologysequence}\dots \to H^q_Y\bigl(Z,\mathscr{O}_Z\bigr) \to H^q \bigl(Z,
\mathscr{O}_Z\bigr) \overset{h_{Z,V}^q}{\longrightarrow} H^q\bigl(V, \mathscr{O}_V\bigr) \to  H^{q+1}_Y\bigl(Z,\mathscr{O}_Z\bigr)\to \cdots.
\end{equation}
\vspace{-5mm}
\begin{equation}\label{difference}
\text{Since $\dim X\hq G \geq 2$, Theorem~\ref{thm:GRvanishingI} yields $H^q_Y\bigl(Z, \mathscr{O}_Z\bigr) = \{0\}$ for $1 \leq q\leq n-1$. \tag{$\star$}}
\end{equation}
Hence, as a consequence of \eqref{localcohomologysequence}, the map $h_{Z,V}^q$ is injective for every $1\leq q \leq n-1$.

The restriction of $p$ to $V = p^{-1}(U)$ is a resolution of $U$. Since $\bigcup_{q\geq 1}\supp R^qp_*\mathscr{O}_Z \subset \{x_0\}$, the complex space $U$ has rational singularities, and the Leray
spectral sequence for $p|_{V}$ yields that $h_{U,V}^q$ is bijective. Similar arguments show that $h_{U',\widetilde U}^q$ is
bijective. Furthermore, Lemma \ref{lem:pullbackinjective} implies that $h_{U, U'}^q$ is injective for all $q \geq 1$.

By the considerations above the map
\[h_{Z,\widetilde{U}}^q \definiere h_{U',\widetilde{U}}^{q} \circ h_{U,U'}^{q} \circ (h^q_{U,V})^{-1} \circ h_{Z,V}^{q}\]
is injective for every $1 \leq q\leq  n-1$. Diagram \eqref{diagramrational} implies $h_{Z, \widetilde{U}}^q = h_{\widetilde{X},
\widetilde{U}}^q \circ h_{Z, \widetilde{X}}^q$, and therefore $h^q_{Z,\widetilde{X}}$ is injective for every $1 \leq q \leq n-1$. Consequently, we
have $H^q\bigl(Z, \mathscr{O}_Z\bigr) = 0$ for all $q\geq 1$. This concludes the proof of part (1) of Theorem \ref{thm:maintheorem}.

The proof of part (2) is completely analogous to the argument given above, except that instead of invoking Theorem~\ref{thm:GRvanishingI} at \eqref{difference} we use Theorem~\ref{thm:GRvanishingII} to show that $h^1_{Z,V}$ is injective.
\end{proof}
Inspecting the proof above gives a refinement of Boutot's result in the algebraic category (see also \cite{GrebSingularities}):
\begin{thm}\label{cor:Boutotrefined}
Let $X$ be a normal algebraic $G$-variety with good quotient $\pi: X \to X\hq G$. Let $f: \widetilde X \to X$ be a resolution of $X$, and let $g: Z \to X\hq G$ a resolution of $X\hq G$. Assume that $R^jf_*\mathscr{O}_{\widetilde X} = 0$ for $1 \leq j \leq q$. Then also $R^jg_*\mathscr{O}_Z = 0$ for $1 \leq j \leq q$.
\end{thm}
\begin{proof}
The proof of Theorem~\ref{thm:maintheorem} (1) given above applies verbatim up to the point~\eqref{difference} where one uses Theorem~\ref{thm:GRvanishingI} to conclude that $H^q_F(Z, \mathscr{O}_Z)$ vanishes (this required the union $\bigcup_{q\geq 1}\supp R^q p_*\mathscr{O}_Z$ to be contained in $\{x_0\}$). However, algebraic and analytic cohomology groups with support always coincide in the algebraic category (\cite[Thm.\ 2.8]{HartshorneLocalCohomology}). Consequently, the cohomology group $H^q_F\bigl(Z, \mathscr{O}_Z\bigr)$ vanishes without the additional assumption on the supports, see \cite{Boutot} and \cite[Prop. 2.2]{HartshorneOgus}.
\end{proof}
This result allows to obtain information about good quotients of algebraic $G$-varieties with worse than rational singularities. See for example \cite{KovacsDuBoisI} for cohomological properties of varieties with Cohen-Macaulay singularities.
\subsection{Quotients of $K$-spaces: proof of Theorem~\ref{thm:mainthmK}}
The proof of Theorem~\ref{thm:mainthmK} follows from the existence of complexifications for Stein $K$-spaces.
\begin{proof}[Proof of Theorem~\ref{thm:mainthmK}]
Let $X$ be a complex $K$-space with analytic Hilbert quotient $\pi: X \to X\hq K$. Since the claim is local on $X\hq K$, we may assume that both $X$ and $X\hq K$ are Stein. Then, by the main result of \cite{HeinznerGIT} there exists a holomorphic Stein $K^\C$-space $X^\C$ and a holomorphic $K$-equivariant open embedding $\imath: X \hookrightarrow X^\C$ with the following properties:
\begin{enumerate}
\item \label{a}the image $\imath(X)$ is an orbit-convex Runge subset of $X^\C$, and $K^\C \acts \imath(X) = X^\C$,

\item \label{b}the analytic Hilbert quotient $\pi^\C: X^\C \to X^\C\negthinspace\hq K^\C$ exists, and the restriction $\pi^\C|_{\imath(X)}: \imath(X) \to X^\C\negthinspace\hq K^\C$ induces an isomorphism $X\hq K \cong X^\C\negthinspace\hq K^\C$ such that the following diagram commutes
    \[\begin{xymatrix}{
    X   \ar^>>>>>>>>\imath[r] \ar_\pi[d]      & X^\C \ar^{\pi^\C}[d]\\
    X\hq K   \ar^>>>>>{\cong}[r] & X^\C\negthinspace\hq K^\C.
    }\end{xymatrix}
    \]
\end{enumerate}
Hence, the claim is proven once we show that the singularities of $X^\C$ are not worse than the singularities of $X$: let $f: \widetilde X \to X^\C$ be a resolution and suppose that $A^q \definiere \supp R^qf_*\mathscr{O}_{\widetilde X}$ is non-empty; since $A^q$ is a $K^\C$-invariant analytic subset of $X^\C$ (cf.\ Section~\ref{sect:sing}) it follows from properties~\eqref{a} and \eqref{b} above that $A^q \cap \imath(X) \neq \emptyset$, contrary to the assumption on the singularities of $X$.
\end{proof}

\end{document}